\documentclass[11pt]{article}

\usepackage{latexsym}
\usepackage{amssymb}

\newtheorem{theorem}{Theorem}
\newtheorem{proposition}{Proposition}
\newtheorem{corollary}{Corollary}

\begin{document}
{
\begin{center}
{\Large\bf
On some hypergeometric Sobolev orthogonal polynomials with several continuous parameters.}
\end{center}
\begin{center}
{\bf S.M. Zagorodnyuk}
\end{center}

\noindent
\textbf{Abstract.}
In this paper we study the following hypergeometric polynomials:
$\mathcal{P}_n(x) = \mathcal{P}_n(x;\alpha,\beta,\delta_1,\dots,\delta_\rho,\kappa_1,\dots,\kappa_\rho) =
{}_{\rho+2} F_{\rho+1} (-n,n+\alpha+\beta+1,\delta_1+1,\dots,\delta_\rho+1;\alpha+1,\kappa_1+\delta_1+1,\dots,\kappa_\rho+\delta_\rho+1;x)$,
and
$\mathcal{L}_n(x) = \mathcal{L}_n(x;\alpha,\delta_1,\dots,\delta_\rho,\kappa_1,\dots,\kappa_\rho) =
{}_{\rho+1} F_{\rho+1} (-n,\delta_1+1,\dots,\delta_\rho+1;\alpha+1,\kappa_1+\delta_1+1,\dots,\kappa_\rho+\delta_\rho+1;x)$,
$n\in\mathbb{Z}_+$, where
$\alpha,\beta,\delta_1,\dots,\delta_\rho\in(-1,+\infty)$, and
$\kappa_1,\dots,\kappa_\rho\in\mathbb{Z}_+$, are some parameters.
The natural number $\rho$ of the continuous parameters $\delta_1,\dots,\delta_\rho$ can be chosen arbitrarily large.
It is seen that the special case $\kappa_1=\dots=\kappa_\rho=0$ leads to Jacobi and Laguerre orthogonal polynomials.
In general, it is shown that polynomials $\mathcal{P}_n(x)$ and $\mathcal{L}_n(x)$ are Sobolev orthogonal polynomials on the real line
with some explicit matrix measures.
We study integral representations, differential equations and generating functions for these polynomials.
Recurrence relations and properties of their zeros are discussed as well.



\section{Introduction.}

The theory of orthogonal polynomials on the real line (OPRL) is a classical subject of analysis having a lot 
of applications~\cite{cit_50000_Gabor_Szego},\cite{cit_50_Bateman},\cite{cit_5000_Ismail}.
The theory of Sobolev orthogonal polynomials is less developed and recognized and it still remains to be a \textit{terra incognita}
in some aspects~\cite{cit_5150_M_X}. 
As this theory may be viewed as a generalization of the classical one, then one can expect that some properties and objects from the classical theory 
will have their mirrors and extensions in the theory of Sobolev orthogonal polynomials. 
For instance, the important property for 
OPRL is that the multiplication by $x$ operator in the corresponding $L^2_\mu$ space is symmetric.
Under some general assumptions, a weaker property of symmetry with respect to an indefinite metric holds
for Sobolev orthogonal polynomials~\cite{cit_90000_Zagorodnyuk_AOT_2022}.
We intend to define and study some generalizations of Jacobi and Laguerre orthogonal polynomials.
Namely, we shall study the following polynomials:
$$ \mathcal{P}_n(x) = \mathcal{P}_n(x;\alpha,\beta,\delta_1,\dots,\delta_\rho,\kappa_1,\dots,\kappa_\rho) = $$
\begin{equation}
\label{f1_5}
= {}_{\rho+2} F_{\rho+1} (-n,n+\alpha+\beta+1,\delta_1+1,\dots,\delta_\rho+1;\alpha+1,\kappa_1+\delta_1+1,\dots,\kappa_\rho+\delta_\rho+1;x),
\end{equation}
and
$$ \mathcal{L}_n(x) = \mathcal{L}_n(x;\alpha,\delta_1,\dots,\delta_\rho,\kappa_1,\dots,\kappa_\rho) = $$
\begin{equation}
\label{f1_7}
= {}_{\rho+1} F_{\rho+1} (-n,\delta_1+1,\dots,\delta_\rho+1;\alpha+1,\kappa_1+\delta_1+1,\dots,\kappa_\rho+\delta_\rho+1;x),\quad
n\in\mathbb{Z}_+, 
\end{equation}
where
$\alpha,\beta,\delta_1,\dots,\delta_\rho\in(-1,+\infty)$, and
$\kappa_1,\dots,\kappa_\rho\in\mathbb{Z}_+$, are some parameters.
Observe that the number $\rho\in\mathbb{N}$ of the continuous parameters $\delta_1,\dots,\delta_\rho$ can be arbitrarily large.
It is clear that the special case $\kappa_1=\dots=\kappa_\rho=0$ leads to the Jacobi and Laguerre orthogonal polynomials on the real line.
There are also some other special cases and related systems of hypergeometric polynomials which were studied before, including Fasenmyer's 
polynomials, see~\cite{cit_5150_Rainville}.
In general, polynomials $\mathcal{P}_n(x)$ and $\mathcal{L}_n(x)$ turns out to be Sobolev orthogonal polynomials on the real line
with some explicit matrix measures. This can be derived on a way proposed in papers~\cite{cit_80000_Zagorodnyuk_JAT_2020}
and~\cite{cit_80000_Zagorodnyuk_JDEA_2021}.

The content of the paper is organized as follows. At first, we shall study generating functions
for the following hypergeometric polynomials:

$$ \mathbf{P}_n(x) = \mathbf{P}_n(x;a,\alpha_1,\dots,\alpha_p;\beta_1,\dots,\beta_q) = $$
\begin{equation}
\label{f1_9}
= {}_{p+2} F_{q} (-n,n+a,\alpha_1,\dots,\alpha_p;\beta_1,\dots,\beta_q;x),
\end{equation}
and
$$ \mathbf{L}_n(x) = \mathbf{L}_n(x;\alpha_1,\dots,\alpha_p;\beta_1,\dots,\beta_q) = $$
\begin{equation}
\label{f1_14}
= {}_{p+1} F_{q} (-n,\alpha_1,\dots,\alpha_p;\beta_1,\dots,\beta_q;x),\quad
n\in\mathbb{Z}_+, 
\end{equation}
where
$a\in(-1,+\infty)$; $\alpha_1,\dots,\alpha_p;\beta_1,\dots,\beta_q\in(0,+\infty)$, are some parameters.
Here $p,q\in\mathbb{Z}_+$, and the case $p=0$ and/or $q=0$ means that $\alpha_k$s and/or $\beta_k$s are absent, respectively.

The generating function for $\mathbf{P}_n(x)$, given by~Theorem~\ref{t2_1} below, seems to be new. It generalizes a (formal) generating function
which appeared in~\cite{cit_70_Bateman} (see formula~(21) on page~266). 
The generating function for $\mathbf{L}_n(x)$, provided by~Theorem~\ref{t2_2}, appeared earlier in~\cite{cit_70_Bateman} (see formula~(25) 
on page~267). 
Notice that it was also a formal one. As far as we know, the convergence questions for the corresponding power series were not studied previously.
The radii which provide the convergence of the corresponding power series will be given by Theorem~\ref{t2_1} (in the case of $\mathbf{P}_n(x)$)
and by Theorem~\ref{t2_2} (in the case of $\mathbf{L}_n(x)$).
In the first case we shall use some series manipulations close to the basic formal series manipulation method, see e.g.~\cite{cit_5700_McBride__1971}.
However our manipulations are not formal, but we take care about the convergence questions. The case of $\mathbf{L}_n$ is simpler, and
we can check the required relation~(\ref{f2_25}) by expanding the generating function into the Taylor series.

Integral representations for the polynomials $\mathcal{P}_n(x)$ and $\mathcal{L}_n(x)$ will be provided by
Corollary~\ref{c2_2} and Proposition~\ref{p2_1}.

Sobolev orthogonality for the polynomials $\mathcal{P}_n(x)$ and $\mathcal{L}_n(x)$ will be obtained in Theorem~\ref{t2_3}.
Here we shall use tools developed earlier in~\cite{cit_80000_Zagorodnyuk_JAT_2020} and~\cite{cit_80000_Zagorodnyuk_JDEA_2021}.
In~\cite{cit_80000_Zagorodnyuk_JDEA_2021} it was shown that Sobolev orthogonal polynomials are related by a differential
equation with orthogonal systems $\mathcal{A}$ of functions acting in the direct sums of usual $L^2_\mu$ spaces 
of square-summable (classes of the equivalence of) functions with respect to a positive measure $\mu$.
The case of a unique $L^2_\mu$ is of a special interest, since it allows to use OPRL to obtain explicit systems of
Sobolev orthogonal polynomials.

The importance of the orthogonality property was our main reason to concentrate our attention on polynomials
$\mathcal{P}_n(x)$ and $\mathcal{L}_n(x)$.
Differential equations for the polynomials $\mathcal{P}_n(x)$ and $\mathcal{L}_n(x)$ will be presented in 
Proposition~\ref{p2_2}.
Known methods for generating functions (see, e.g., \cite[Chapter XIX]{cit_70_Bateman}, \cite{cit_5700_McBride__1971}) can be used 
to obtain some additional properties of the polynomials $\mathcal{P}_n(x)$ and $\mathcal{L}_n(x)$.
We shall discuss the existence of recurrence relations for these polynomials.
In Theorem~\ref{t2_4} we obtain a five-term recurrence relation for a special case of polynomials $\mathbf{L}_n(x)$,
with $p=2$, $q=3$. The latter provides a five-term recurrence relation for $\mathcal{L}_n$ with $\rho=2$, as a special case.
In this case the polynomials $\mathcal{L}_n(x)$ ($\rho =2$) have three important properties: 
\begin{itemize}
\item[(1)] the Sobolev orthogonality;

\item[(2)] these polynomials are 
(generalized) eigenvalues of a pencil of differential operators;

\item[(3)] these polynomials are eigenvalues of a pencil of difference operators.
\end{itemize}

Of course, each of these features is valuable and $\mathcal{L}_n$ ($\rho =2$) possess all of them. 
These properties make polynomials $\mathcal{L}_n(x)$ close to classical systems of polynomials
and their generalizations, see \cite{cit_50000_Gabor_Szego},\cite{cit_5105_Koekoek_book}.
Observe that properties $(2)$ and $(3)$ are close to the bispectral problems studied
for various orthogonal systems of functions, 
see~\cite{cit_DG_1986},\cite{cit_110_Everitt_2001},\cite{cit_7100__S_Zh},\cite{cit210_Horozov},\cite{cit_110_Duran} 
and references therein.

Finally, some information on the location of zeros for $\mathcal{P}_n(x)$ and $\mathcal{L}_n(x)$ will be given in Proposition~\ref{p2_3}.

\noindent
{\bf Notations. }
As usual, we denote by $\mathbb{R}, \mathbb{C}, \mathbb{N}, \mathbb{Z}, \mathbb{Z}_+$,
the sets of real numbers, complex numbers, positive integers, integers and non-negative integers,
respectively; $\mathbb{D}_r := \{ z\in\mathbb{C}:\ |z|<r \}$, $r>0$; $\mathbb{D} := \mathbb{D}_1$. 
By $\mathbb{Z}_{k,l}$ we mean all integers $j$ satisfying the following inequality:
$k\leq j\leq l$; ($k,l\in\mathbb{Z}$).
By $\mathbb{P}$ we denote the set of all polynomials with complex coefficients.
By $\mathbb{P}_r$ we mean the set of all polynomials with real coefficients.
By $M^T$ we mean the transpose of a complex matrix $M$.
For a complex number $c$ we denote
$(c)_0 = 1$, $(c)_k = c\cdots (c+k-1)$, $k\in\mathbb{N}$ (\textit{the shifted factorial or Pochhammer's symbol}).
As usual, the generalized hypergeometric function is denoted by
$$ {}_m F_n(a_1,\dots,a_m; b_1,\dots,b_n;x) = 
{}_m F_n \left[
\begin{array}{cc} a_1,\dots,a_m; \\
b_1,\dots,b_n;
\end{array}
x
\right] = $$
$$ = \sum_{k=0}^\infty \frac{(a_1)_k ... (a_m)_k}{(b_1)_k ... (b_n)_k} \frac{x^k}{k!}, $$
where $a_j$, $b_j$, $x$ are complex numbers and $b_j$s are not allowed to take negative integer values.
By $\Gamma(z)$ and $\mathrm{B}(z)$ we denote the gamma function and the beta function, respectively.

\section{Properties of some hypergeometric Sobolev orthogonal polynomials.}

The following theorem provides a generating function for the polynomials $\mathbf{P}_n(x)$ (cf.~\cite[p. 266]{cit_70_Bateman}, formula~(21)).

\begin{theorem}
\label{t2_1}
Let $p,q\in\mathbb{Z}_+$: $p\leq q-1$, and  $c$: $0<c<\frac{1}{2}$, be fixed.
Let $a;\alpha_1,\dots,\alpha_p;\beta_1,\dots,\beta_q\in(0,+\infty)$, be arbitrary parameters.
The following relation holds:
$$ (1-t)^{-a} 
{}_{p+2} F_{q} \left( 
\frac{a}{2}, \frac{a+1}{2},\alpha_1,\dots,\alpha_p;\beta_1,\dots,\beta_q; -\frac{4xt}{ (1-t)^2 }  
\right)
= $$
\begin{equation}
\label{f2_5}
= \sum_{n=0}^\infty \frac{ (a)_n }{ n! } \mathbf{P}_n(x;a,\alpha_1,\dots,\alpha_p;\beta_1,\dots,\beta_q) t^n,
\end{equation}
where
\begin{equation}
\label{f2_7}
t,x\in\mathbb{C}:\ |t|<c,\ |x|< \frac{1}{4c} - \frac{1}{2}.
\end{equation}
 
\end{theorem}

\noindent
\textbf{Proof.} 
Notice that condition~(\ref{f2_7}) provides that
\begin{equation}
\label{f2_9}
\left|
\frac{4xt}{ (1-t)^2 }
\right| < 1.
\end{equation}
In fact, we may write:
$$ \left|
\frac{4xt}{ (1-t)^2 }
\right| =
\frac{ 4 |x| |t| }{ |1-t|^2 } < 
\frac{ 4 \left( \frac{1}{4c} - \frac{1}{2} \right) c }{ (1-c)^2 }
= \frac{1-2c}{(1-c)^2} \leq \frac{1-2c + c^2 }{(1-c)^2} = 1. $$
Therefore the left-hand side of~(\ref{f2_5}) is well-defined for all $t,x$ satisfying condition~(\ref{f2_7}).
Denote by $R_1$ the right-hand side of~(\ref{f2_5}). At this point we do not know if the series in $R_1$ converges. 
Consider the following two iterated series which differ by the order of summation:
\begin{equation}
\label{f2_14}
R_2 := \sum_{n=0}^\infty \sum_{k=0}^\infty (a)_n \frac{t^n}{n!} (-n)_k (n+a)_k u_k \frac{x^k}{k!},  
\end{equation}
\begin{equation}
\label{f2_16}
R_3 := \sum_{k=0}^\infty \sum_{n=0}^\infty (a)_n \frac{t^n}{n!} (-n)_k (n+a)_k u_k \frac{x^k}{k!},  
\end{equation}
where for brevity we denoted
\begin{equation}
\label{f2_17}
u_j := \frac{ (\alpha_1)_j \dots (\alpha_p)_j }{ (\beta_1)_j \dots (\beta_q)_j },\qquad j\in\mathbb{Z}_+,
\end{equation}
and $t,x$ are satisfying condition~(\ref{f2_7}).
We are going to prove that the series $R_3$ converges absolutely. 
Then we could conclude that the series $R_2$
converges and the sums of $R_2$ and $R_3$ coincide (see, e.g., Theorem~3 in~\cite[p. 34]{cit_90_Fichtenholz}). After this will be done, we shall check 
that $R_3$ coincides with the left-hand side of~(\ref{f2_5}). On the other hand, relation $R_2=R_1$ will be obviously satisfied.
Denote
$$ \widehat R_3 := \sum_{k=0}^\infty \sum_{n=0}^\infty \left| (a)_n \frac{t^n}{n!} (-n)_k (n+a)_k u_k \frac{x^k}{k!} \right| = $$
$$ = \sum_{k=0}^\infty u_k \frac{|x|^k}{k!} \sum_{n=k}^\infty (a)_n (n+a)_k  |(-n)_k|  \frac{|t|^n}{n!}  = $$
\begin{equation}
\label{f2_20}
= \sum_{k=0}^\infty u_k \frac{|x|^k}{k!} \sum_{n=k}^\infty (a)_n (n+a)_k  \frac{ |t|^n }{ (n-k)! }, 
\end{equation}
where we have removed the null terms.
Denote the inner sum in the last row of~(\ref{f2_20}) by $S_k$.
By the ratio test it converges for all $t\in\mathbb{D}_c$.
Changing the summation index $j=n-k$ we get
$$ S_k = \sum_{j=0}^\infty (a)_{j+k} (j+k+a)_k  \frac{ |t|^{j+k} }{ j! } = 
(a)_{2k} |t|^k \sum_{j=0}^\infty (a+2k)_j \frac{|t|^j}{j!} =  $$ 
$$ = (a)_{2k}|t|^k (1-|t|)^{-a-2k},\qquad  t\in\mathbb{D}_c. $$
Then
$$ \widehat R_3 = (1-|t|)^{-a} \sum_{k=0}^\infty (a)_{2k} \frac{u_k}{k!} \left( \frac{|xt|}{ (1-|t|)^2 }  \right)^k  = $$
$$ = (1-|t|)^{-a} \sum_{k=0}^\infty \left( \frac{a}{2} \right)_k \left( \frac{a+1}{2} \right)_k \frac{u_k}{k!} \left( \frac{4|xt|}{ (1-|t|)^2 }  
\right)^k  = $$
\begin{equation}
\label{f2_22}
= (1-|t|)^{-a} 
{}_{p+2} F_{q} \left( 
\frac{a}{2}, \frac{a+1}{2},\alpha_1,\dots,\alpha_p;\beta_1,\dots,\beta_q; \frac{4|x||t|}{ (1-|t|)^2 }  
\right),
\end{equation}
where we have used the following relation (see~Lemma~5 in~\cite[p. 22]{cit_5150_Rainville}):
$$ (a)_{2k} = 4^k \left( \frac{a}{2} \right)_k \left( \frac{a+1}{2} \right)_k. $$
By virtue of~(\ref{f2_9}) with parameters $|x|,|t|$ instead of $x,t$, we obtain that
$\frac{ 4|x||t| }{ (1-|t|)^2 } < 1$, and this proves the last line of~(\ref{f2_22}).
Thus, the series $R_3$ converges absolutely. Let
$$ R_3 = \sum_{k=0}^\infty \sum_{n=0}^\infty a_{k,n},\quad a_{k,n} = u_{k,n} + i v_{k,n},\ u_{k,n}, v_{k,n}\in\mathbb{R}. $$
By Theorem~2 in~\cite[p. 34]{cit_90_Fichtenholz} we conclude that
$$ \sum_{j=0}^\infty | a_j | < \infty, $$
where the series is composed of elements $a_{k,j}$, placed in an arbitrary order.
Let $a_j = u_j + i v_j$, $u_j, v_j\in\mathbb{R}$.
By the comparison test it follows that
$$ \sum_{j=0}^\infty | u_j | < \infty,\ \sum_{j=0}^\infty | v_j | < \infty. $$
By Theorem~1 in~\cite[p. 32]{cit_90_Fichtenholz} we obtain that
\begin{equation}
\label{f2_22_1}
\sum_{k=0}^\infty \sum_{n=0}^\infty u_{k,n} = \sum_{n=0}^\infty \sum_{k=0}^\infty u_{k,n} = \sum_{j=0}^\infty u_j;
\end{equation}
\begin{equation}
\label{f2_22_2}
\sum_{k=0}^\infty \sum_{n=0}^\infty i v_{k,n} = \sum_{n=0}^\infty \sum_{k=0}^\infty i v_{k,n} = \sum_{j=0}^\infty i v_j.
\end{equation}
Summing relations~(\ref{f2_22_1}) and~(\ref{f2_22_2}) we get
\begin{equation}
\label{f2_22_3}
\sum_{k=0}^\infty \sum_{n=0}^\infty a_{k,n} = \sum_{n=0}^\infty \sum_{k=0}^\infty a_{k,n}.
\end{equation}
Therefore $R_3 = R_2$.
It remains to check that $R_3$ coincides with the left-hand side of~(\ref{f2_5}).
We may write:
$$ R_3 =  \sum_{k=0}^\infty  u_k \frac{x^k}{k!} \sum_{n=k}^\infty (a)_n (-n)_k (n+a)_k \frac{t^n}{n!}.  $$
Denote 
$$ T_k := \sum_{n=k}^\infty (a)_n (-n)_k (n+a)_k \frac{t^n}{n!}. $$
The series $T_k$ converges absolutely by the ratio test.
Proceeding in a similar manner as for $S_k$, we change the summation index $j=n-k$:
$$ T_k = \sum_{j=0}^\infty (a)_{j+k} (j+k+a)_k  (-1)^k \frac{ t^{j+k} }{ j! } = 
(a)_{2k} (-t)^k \sum_{j=0}^\infty (a+2k)_j \frac{ t^j }{j!} =  $$ 
$$ = (a)_{2k} (-t)^k (1-t)^{-a-2k},\qquad  t\in\mathbb{D}_c. $$
Therefore
$$ R_3 = (1-t)^{-a} \sum_{k=0}^\infty (a)_{2k} \frac{u_k}{k!} \left( \frac{ -xt }{ (1-t)^2 }  \right)^k  = $$
$$ = (1-t)^{-a} \sum_{k=0}^\infty \left( \frac{a}{2} \right)_k \left( \frac{a+1}{2} \right)_k \frac{u_k}{k!} \left( -\frac{4xt}{ (1-t)^2 }  
\right)^k  = $$
\begin{equation}
\label{f2_24}
= (1-t)^{-a} 
{}_{p+2} F_{q} \left( 
\frac{a}{2}, \frac{a+1}{2},\alpha_1,\dots,\alpha_p;\beta_1,\dots,\beta_q; \frac{ -4xt }{ (1-t)^2 }  
\right),
\end{equation}
where we have used relation~(\ref{f2_9}).
Since $R_3 = R_2 = R_1$, the proof is complete.
$\Box$

Now we shall obtain a generating function for the polynomials $\mathbf{L}_n(x)$ (cf.~\cite[p. 267]{cit_70_Bateman}, formula~(25)).

\begin{theorem}
\label{t2_2}
Let $p,q\in\mathbb{Z}_+$: $p\leq q+1$, be fixed.
Let $\alpha_1,\dots,\alpha_p;\beta_1,\dots,\beta_q\in(0,+\infty)$, be arbitrary parameters.
The following relation holds:
\begin{equation}
\label{f2_25}
e^t 
{}_{p} F_{q} \left( 
\alpha_1,\dots,\alpha_p;\beta_1,\dots,\beta_q; -xt
\right) = \sum_{n=0}^\infty \mathbf{L}_n(x;\alpha_1,\dots,\alpha_p;\beta_1,\dots,\beta_q) \frac{t^n}{ n! },
\end{equation}
where $t,x\in\mathbb{D}$.
If $p\leq q$ then relation~(\ref{f2_25}) holds for all $t,x\in\mathbb{C}$.
\end{theorem}

\noindent
\textbf{Proof.} 
Denote by $g(t) = g_x(t)$ the left-hand side of~(\ref{f2_25}).
Set 
$$ D := \left\{
\begin{array}{cc} \mathbb{D}, & \mbox{if } p = q+1 \\
\mathbb{C}, & \mbox{if } p \leq q \end{array} 
\right.. $$
Fix an arbitrary $x\in D$. Then $g(t)=g_x(t)$ is an analytic function of $t$ in the domain $D$.
Let us calculate Taylor's coefficients for its expansion at $t=0$. By the Leibniz rule we may write:
$$ g^{(n)} (0) = \sum_{k=0}^n {n \choose k}
\left.
\left(
{}_{p} F_{q} \left( 
\alpha_1,\dots,\alpha_p;\beta_1,\dots,\beta_q; -xt
\right)
\right)^{(k)}_t \right|_{t=0} 
\left.
\left(
e^t
\right)^{(n-k)} \right|_{t=0} = $$
$$ = \sum_{k=0}^n {n \choose k}
\left.
\left(
\sum_{j=0}^\infty \frac{ (\alpha_1)_j \dots (\alpha_p)_j  }{ (\beta_1)_j \dots (\beta_q)_j }
\frac{ (-x)^j }{ j! } t^j
\right)^{(k)}_t \right|_{t=0} = $$
$$ = \sum_{k=0}^n {n \choose k}
\frac{ (\alpha_1)_k \dots (\alpha_p)_k  }{ (\beta_1)_k \dots (\beta_q)_k }
(-x)^k
= \sum_{k=0}^n (-n)_k
\frac{ (\alpha_1)_k \dots (\alpha_p)_k  }{ (\beta_1)_k \dots (\beta_q)_k }
\frac{x^k}{k!}
= $$
$$ = 
\mathbf{L}_n(x;\alpha_1,\dots,\alpha_p;\beta_1,\dots,\beta_q). $$
Thus, relation~(\ref{f2_25}) coincides with Taylor's expansion of $g(t)$ at $t=0$.
$\Box$

Observe that
$$ \mathcal{P}_n(x;\alpha,\beta,\delta_1,\dots,\delta_\rho,\kappa_1,\dots,\kappa_\rho) = $$
\begin{equation}
\label{f2_35}
= \mathbf{P}_n(x;\alpha+\beta+1,\delta_1 +1,\dots,\delta_\rho +1;\alpha+1,\kappa_1+\delta_1 +1,\dots,\kappa_\rho+\delta_\rho +1),\quad
n\in\mathbb{Z}_+, 
\end{equation}
and
$$ \mathcal{L}_n(x;\alpha,\delta_1,\dots,\delta_\rho,\kappa_1,\dots,\kappa_\rho) = $$
\begin{equation}
\label{f2_37}
= \mathbf{L}_n(x;\delta_1 +1,\dots,\delta_\rho +1;\alpha+1,\kappa_1+\delta_1 +1,\dots,\kappa_\rho+\delta_\rho +1),\quad
n\in\mathbb{Z}_+, 
\end{equation}
where
$\alpha,\beta,\delta_1,\dots,\delta_\rho\in(-1,+\infty)$, and
$\kappa_1,\dots,\kappa_\rho\in\mathbb{Z}_+$, are arbitrary parameters; $\rho\in\mathbb{N}$.
Thus, we can formulate the following corollary.

\begin{corollary}
\label{c2_1}
Let $\rho\in\mathbb{N}$, and 
$\delta_1,\dots,\delta_\rho\in(-1,+\infty)$; $\kappa_1,\dots,\kappa_\rho\in\mathbb{Z}_+$, be arbitrary parameters.
The following statements hold:

\noindent
$(i)$ Let $c:\ 0<c<\frac{1}{2}$, and $\alpha,\beta\in(-1,+\infty):\ \alpha + \beta > -1$, be given. Then
$$ (1-t)^{-\alpha-\beta-1} 
{}_{\rho+2} F_{\rho+1} 
\left[
\begin{array}{cc} \frac{\alpha+\beta+1}{2}, \frac{\alpha+\beta+2}{2},\delta_1+1,\dots,\delta_\rho+1; \\
\alpha+1,\kappa_1+\delta_1+1,\dots,\kappa_\rho+\delta_\rho+1;
\end{array}
-\frac{4xt}{ (1-t)^2 }  
\right]
= $$
\begin{equation}
\label{f2_40}
= \sum_{n=0}^\infty \frac{ (\alpha+\beta+1)_n }{ n! } \mathcal{P}_n(x;\alpha,\beta,\delta_1,\dots,\delta_\rho,\kappa_1,\dots,\kappa_\rho) t^n,
\end{equation}
where
\begin{equation}
\label{f2_46}
t,x\in\mathbb{C}:\ |t|<c,\ |x|< \frac{1}{4c} - \frac{1}{2};
\end{equation}

\noindent
$(ii)$ Let $\alpha > -1$. For all $t,x\in\mathbb{C}$ the following relation is valid:
$$ e^t 
{}_{\rho} F_{\rho+1} 
\left[
\begin{array}{cc} \delta_1+1,\dots,\delta_\rho+1; \\
\alpha+1,\kappa_1+\delta_1+1,\dots,\kappa_\rho+\delta_\rho+1;
\end{array}
-xt
\right]
 = $$
\begin{equation}
\label{f2_48}
= \sum_{n=0}^\infty \mathcal{L}_n(x;\alpha,\delta_1,\dots,\delta_\rho,\kappa_1,\dots,\kappa_\rho) \frac{t^n}{ n! }.
\end{equation}

\end{corollary}

\noindent
\textbf{Proof.} 
By Theorems~\ref{t2_1} and~\ref{t2_2} the proof is straightforward.
$\Box$

\begin{corollary}
\label{c2_2}
Let $\rho\in\mathbb{N}$, and 
$\delta_1,\dots,\delta_\rho\in(-1,+\infty)$; $\kappa_1,\dots,\kappa_\rho\in\mathbb{Z}_+$, be arbitrary parameters.
If $\alpha>-1$, then 
$$ \mathcal{L}_n(x;\alpha,\delta_1,\dots,\delta_\rho,\kappa_1,\dots,\kappa_\rho) = $$
$$ = \frac{ n! }{ 2\pi i } \oint_{ |\zeta|=1 }
\zeta^{-n-1} e^\zeta
{}_{\rho} F_{\rho+1} 
\left[
\begin{array}{cc} \delta_1+1,\dots,\delta_\rho+1; \\
\alpha+1,\kappa_1+\delta_1+1,\dots,\kappa_\rho+\delta_\rho+1;
\end{array}
-x\zeta
\right] d\zeta,
$$
\begin{equation}
\label{f2_48_3}
x\in\mathbb{C},\quad n\in\mathbb{Z}_+. 
\end{equation}
If $\alpha,\beta\in (-1,+\infty):\ \alpha + \beta > -1$, then
$$ \mathcal{P}_n(x;\alpha,\beta,\delta_1,\dots,\delta_\rho,\kappa_1,\dots,\kappa_\rho) = $$
$$ = \frac{1}{ 2\pi i } \frac{ n! }{ (\alpha+\beta+1)_n } \oint_{ |\zeta|=\frac{1}{4} }
\zeta^{-n-1} (1-\zeta)^{-\alpha-\beta-1} * $$
$$ * {}_{\rho+2} F_{\rho+1} 
\left[
\begin{array}{cc} \frac{\alpha+\beta+1}{2}, \frac{\alpha+\beta+2}{2},\delta_1+1,\dots,\delta_\rho+1; \\
\alpha+1,\kappa_1+\delta_1+1,\dots,\kappa_\rho+\delta_\rho+1;
\end{array}
-\frac{4x\zeta}{ (1-\zeta)^2 }  
\right]
d\zeta, $$
\begin{equation}
\label{f2_48_1}
x\in\mathbb{C}:\ |x| < \frac{1}{4},\quad n\in\mathbb{Z}_+.
\end{equation}

\end{corollary}

\noindent
\textbf{Proof.} 
The proof follows from Corollary~\ref{c2_1}, if one calculate the corresponding Taylor coefficients (with $c=\frac{1}{3}$).
$\Box$

Polynomials $\mathcal{P}_n$ and $\mathcal{L}_n$ also admit some recursive integral representations.
Let $\alpha,\beta>-1$. Consider the classical Jacobi and Laguerre polynomials:
\begin{equation}
\label{f2_48_5}
J_n(x) = J_n(x;\alpha,\beta) :=  {}_2 F_1 (-n,n+\alpha+\beta+1;\alpha+1;x),
\end{equation}
\begin{equation}
\label{f2_48_6}
L_n(x) = L_n(x;\alpha) :=  {}_1 F_1 (-n;\alpha+1;x),\qquad  n\in\mathbb{Z}_+.
\end{equation}

\begin{proposition}
\label{p2_1}
Let $\rho\in\mathbb{N}$, and 
$\alpha,\beta, \delta_1,\dots,\delta_\rho\in(-1,+\infty)$; $\kappa_1,\dots,\kappa_\rho\in\mathbb{N}$, be arbitrary parameters.
If $\rho\geq 2$, then
$$ \mathcal{P}_n(z;\alpha,\beta,\delta_1,\dots,\delta_\rho,\kappa_1,\dots,\kappa_\rho) = $$
$$ = \frac{ \Gamma(\kappa_\rho+\delta_\rho+1) }{ \Gamma(\delta_\rho+1) \Gamma(\kappa_\rho) }
\int_0^1 t^{\delta_\rho} (1-t)^{\kappa_\rho-1}
\mathcal{P}_n(zt;\alpha,\beta,\delta_1,\dots,\delta_{\rho-1},\kappa_1,\dots,\kappa_{\rho-1})dt,
$$
\begin{equation}
\label{f2_48_7}
z\in\mathbb{C}:\ |z|<1,\quad n\in\mathbb{Z}_+. 
\end{equation}
If $\rho= 1$, then
$$ \mathcal{P}_n(z;\alpha,\beta,\delta_1,\kappa_1) = 
\frac{ \Gamma(\kappa_1+\delta_1+1) }{ \Gamma(\delta_1+1) \Gamma(\kappa_1) }
\int_0^1 t^{\delta_1} (1-t)^{\kappa_1-1}
J_n(zt;\alpha,\beta) dt,
$$
\begin{equation}
\label{f2_48_9}
z\in\mathbb{C}:\ |z|<1,\quad n\in\mathbb{Z}_+. 
\end{equation}
If $\rho\geq 2$, then
$$ \mathcal{L}_n(z;\alpha,\delta_1,\dots,\delta_\rho,\kappa_1,\dots,\kappa_\rho) = $$
$$ = \frac{ \Gamma(\kappa_\rho+\delta_\rho+1) }{ \Gamma(\delta_\rho+1) \Gamma(\kappa_\rho) }
\int_0^1 t^{\delta_\rho} (1-t)^{\kappa_\rho-1}
\mathcal{L}_n(zt;\alpha,\delta_1,\dots,\delta_{\rho-1},\kappa_1,\dots,\kappa_{\rho-1})dt,
$$
\begin{equation}
\label{f2_48_11}
z\in\mathbb{C},\quad n\in\mathbb{Z}_+. 
\end{equation}
If $\rho= 1$, then
$$ \mathcal{L}_n(z;\alpha,\delta_1,\kappa_1) = 
\frac{ \Gamma(\kappa_1+\delta_1+1) }{ \Gamma(\delta_1+1) \Gamma(\kappa_1) }
\int_0^1 t^{\delta_1} (1-t)^{\kappa_1-1}
L_n(zt;\alpha) dt,
$$
\begin{equation}
\label{f2_48_14}
z\in\mathbb{C},\quad n\in\mathbb{Z}_+. 
\end{equation}

\end{proposition}

\noindent
\textbf{Proof.} Use hypergeometric representations of the corresponding polynomials and
Theorem~28 in~\cite[p. 85]{cit_5150_Rainville}. 
$\Box$

Fix an arbitrary $\rho\in\mathbb{N}$, and choose arbitrary parameters 
$\alpha,\beta,\delta_1,\dots,\delta_\rho\in(-1,+\infty)$, and
$\kappa_1,\dots,\kappa_\rho\in\mathbb{N}$.
Introduce the following linear differential operator $L = L(\delta,k)$ with 
polynomial coefficients, $\delta>-1$, $k\in\mathbb{N}$:
\begin{equation}
\label{f2_50}
L y(x) = \frac{1}{ (\delta+1) \dots (\delta + k) }
x^{-\delta}
\left(
x^{k + \delta} y(x)
\right)^{(k)},\qquad  y(x)\in\mathbb{P}.
\end{equation}
Denote
$$ \widehat D = \widehat D (\delta_1,\dots,\delta_\rho;\kappa_1,\dots,\kappa_\rho) =
L(\delta_1,\kappa_1)
L(\delta_2,\kappa_2) 
\dots
L(\delta_\rho,\kappa_\rho) =
$$
\begin{equation}
\label{f2_52}
= \sum_{j=0}^\kappa c_j(x) \frac{ d^j }{ dx^j },\qquad  c_j(x)=c_j(x;\delta_1,\dots,\delta_\rho;\kappa_1,\dots,\kappa_\rho)\in\mathbb{P}, 
\end{equation}
where
$c_\kappa(x)$ is not the null polynomial,
$\kappa := \kappa_1 + \dots +\kappa_\rho$.

Now we shall show that the polynomials $\mathcal{P}_n(x)$ and $\mathcal{L}_n(x)$ are Sobolev orthogonal polynomials on the real line.

\begin{theorem}
\label{t2_3}
Let $\rho\in\mathbb{N}$, and
$\alpha,\beta,\delta_1,\dots,\delta_\rho\in(-1,+\infty)$;
$\kappa_1,\dots,\kappa_\rho\in\mathbb{N}$ be arbitrary parameters.
Let $\widehat D = \widehat D (\delta_1,\dots,\delta_\rho;\kappa_1,\dots,\kappa_\rho)$ be given
by~(\ref{f2_52}), and
$$ M(x) := (c_0(x),\dots,c_\kappa(x))^T (c_0(x),\dots,c_\kappa(x)),\qquad x\in\mathbb{R}. $$
For polynomials $\mathcal{P}_n(x)$ and $\mathcal{L}_n(x)$, defined as in~(\ref{f1_5}),(\ref{f1_7}), the following relations hold:
$$ \int_{0}^1 (\mathcal{P}_n(x), \mathcal{P}_n'(x), \dots, \mathcal{P}_n^{(\kappa)}(x)) M(x)
\left(
\begin{array}{cccc}  \mathcal{P}_m(x) \\
\mathcal{P}_m'(x) \\
\vdots \\
\mathcal{P}_m^{(\kappa)}(x) \end{array}
\right)
(1-x)^\alpha
(1+x)^\beta
dx = $$
\begin{equation}
\label{f2_57}
= A_n\delta_{n,m},\qquad  A_n>0,\ n,m\in\mathbb{Z}_+;
\end{equation}

$$ \int_{0}^\infty (\mathcal{L}_n(x), \mathcal{L}_n'(x), \dots, \mathcal{L}_n^{(\kappa)}(x)) M(x)
\left(
\begin{array}{cccc}  \mathcal{L}_m(x) \\
\mathcal{L}_m'(x) \\
\vdots \\
\mathcal{L}_m^{(\kappa)}(x) \end{array}
\right)
x^\alpha e^{-x}
dx = $$
\begin{equation}
\label{f2_59}
= B_n\delta_{n,m},\qquad  B_n>0,\ n,m\in\mathbb{Z}_+.
\end{equation}

\end{theorem}

\noindent
\textbf{Proof.} A direct calculation shows that
$$ L(\delta_\rho,\kappa_\rho) \mathcal{P}_n(x;\alpha,\beta,\delta_1,\dots,\delta_\rho,\kappa_1,\dots,\kappa_\rho) = $$
$$ = \left\{
\begin{array}{cc} 
\mathcal{P}_n(x;\alpha,\beta,\delta_1,\dots,\delta_{\rho-1},\kappa_1,\dots,\kappa_{\rho-1} ), & \mbox{if $\rho\geq 2$} \\
{}_2 F_1 (-n,n+\alpha+\beta+1;\alpha+1;x), & \mbox{if $\rho = 1$}
\end{array}
\right.;
$$
and
$$ L(\delta_\rho,\kappa_\rho) \mathcal{L}_n(x;\alpha,\delta_1,\dots,\delta_\rho,\kappa_1,\dots,\kappa_\rho) = $$
$$ = \left\{
\begin{array}{cc} 
\mathcal{L}_n(x;\alpha,\delta_1,\dots,\delta_{\rho-1},\kappa_1,\dots,\kappa_{\rho-1} ), & \mbox{if $\rho\geq 2$} \\
{}_1 F_1 (-n;\alpha+1;x), & \mbox{if $\rho = 1$}
\end{array}
\right..
$$
Therefore
$$ \widehat D \mathcal{P}_n(x;\alpha,\beta,\delta_1,\dots,\delta_\rho,\kappa_1,\dots,\kappa_\rho) =
{}_2 F_1 (-n,n+\alpha+\beta+1;\alpha+1;x) = J_n(x;\alpha,\beta), $$
and
$$ \widehat D \mathcal{L}_n(x;\alpha,\delta_1,\dots,\delta_\rho,\kappa_1,\dots,\kappa_\rho) =
{}_1 F_1 (-n;\alpha+1;x) = L_n(x;\alpha). $$
The latter expressions for the Jacobi polynomials $J_n$ and the Laguerre polynomials $L_n$
can be inserted into their orthogonality relations to obtain relations~(\ref{f2_57}),(\ref{f2_59}).
This finishes the proof. $\Box$

Of course, the hypergeometric nature of polynomials $\mathcal{P}_n$ and $\mathcal{L}_n$ provides differential
equations for them.

\begin{proposition}
\label{p2_2}
Let $\rho\in\mathbb{N}$, and 
$\alpha,\beta, \delta_1,\dots,\delta_\rho\in(-1,+\infty)$; $\kappa_1,\dots,\kappa_\rho\in\mathbb{Z}_+$, be arbitrary parameters.
Let $\theta = z \frac{d}{dz}$, and
\begin{equation}
\label{f2_65}
K := \theta (\theta + \alpha) \prod_{j=1}^\rho (\theta +\kappa_j + \delta_j),\quad
L := \prod_{k=1}^\rho (\theta +\delta_k + 1),
\end{equation}

\begin{equation}
\label{f2_67}
D_0 := K - z\theta (\theta + \alpha + \beta + 1) L,\quad
D_1 := z L,\quad D_2 := K - z\theta L.
\end{equation}
Then $\forall n\in\mathbb{Z}_+$, 
\begin{equation}
\label{f2_68}
D_0 \mathcal{P}_n(z) = -n(n+\alpha +\beta +1) D_1 \mathcal{P}_n(z),\quad z\in\mathbb{D};
\end{equation}

\begin{equation}
\label{f2_73}
D_2 \mathcal{L}_n(z) = -n D_1 \mathcal{L}_n(z),\quad z\in\mathbb{C}. 
\end{equation}

\end{proposition}

\noindent
\textbf{Proof.} Use hypergeometric representations of the corresponding polynomials and
the differential equation for ${}_p F_q$.
$\Box$

Let us turn to the question of the existence of some recurrence relations for polynomials $\mathbf{P}_n$ and $\mathbf{L}_n$.
Here we shall use once more the powerful tool of generating functions.
We can use the differential equations for $\mathbf{P}_n$ and $\mathbf{L}_n$. However, in the case
of $\mathbf{P}_n$ we have a generating function (see formula~(\ref{f2_5})) which involves
a hypergeometric function with an argument $-\frac{4xt}{ (1-t)^2 }$. For big values of $p$ and $q$ this causes very complicated expressions
during the differentiating, if we use, for example, Fa\`{a} di Bruno's formula.

In the case of $\mathbf{L}_n$ it appears a similar problem. For big values of $p$ and $q$ the expressions for the coefficients
of recurrence relations will be complicated and it is not clear that they will be nontrivial. Thus, the non-triviality
of the recurrence relations can not be guaranteed.

We are not ready to treat effectively the case of general $p$ and $q$. It looks reasonable to
investigate concrete systems of polynomials $\mathbf{P}_n$ or $\mathbf{L}_n$, having some fixed values of $p$ and $q$.
Even in this case expressions for the coefficients can be huge and probably of few use.
We shall study the case $p=2$, $q=3$, for the polynomials $\mathbf{L}_n$:

\begin{equation}
\label{f2_75}
\mathbf{L}_n(x) = \mathbf{L}_n(x;\alpha_1,\alpha_2;\beta_1,\beta_2,\beta_3) 
= {}_{3} F_{3} (-n,\alpha_1,\alpha_2;\beta_1,\beta_2,\beta_3;x),\quad
n\in\mathbb{Z}_+, 
\end{equation}
where
$\alpha_1,\alpha_2,\beta_1,\beta_2,\beta_3\in(0,+\infty)$.
By Theorem~\ref{t2_2} we may write:
\begin{equation}
\label{f2_77}
e^t 
{}_{2} F_{3} \left( 
\alpha_1,\alpha_2;\beta_1,\beta_2,\beta_3; -xt
\right) = 
\sum_{n=0}^\infty \mathbf{L}_n(x) \frac{t^n}{ n! },\quad t,x\in\mathbb{C}.
\end{equation}

Fix an arbitrary number $x\in\mathbb{C}\backslash\{ 0 \}$. Introduce a new variable $z$:
$$ z = -x t. $$
Relation~(\ref{f2_77}) may be written in the following form:
\begin{equation}
\label{f2_78}
{}_{2} F_{3} \left( 
\alpha_1,\alpha_2;\beta_1,\beta_2,\beta_3; z
\right) = 
e^{ \frac{z}{x} }
\sum_{n=0}^\infty \mathbf{L}_n(x) \frac{ (-1)^n }{ x^n }\frac{z^n}{ n! },\quad z\in\mathbb{C}.
\end{equation}

Denote the left-hand side of relation~(\ref{f2_78}) by $w(z)$. It satisfies the differential equation for the hypergeometric
function:
\begin{equation}
\label{f2_79}
\left[
\theta
(\theta +\beta_1 -1) (\theta +\beta_2 -1) (\theta +\beta_3 -1)
-z 
(\theta + \alpha_1)(\theta + \alpha_2)
\right] w(z) = 0,
\end{equation}
where $\theta = z\frac{d}{dz}$.
Set
\begin{equation}
\label{f2_82}
b_1 := \beta_1 - 1,\quad b_2 := \beta_2 - 1,\quad b_3 := \beta_3 - 1,  
\end{equation}
\begin{equation}
\label{f2_84}
c := b_1 + b_2 + b_3 + 6,\quad 
\widehat b := 7 + 3 ( b_1 + b_2 + b_3 ) + b_1 b_2 + b_1 b_3 + b_2 b_3,
\end{equation}
\begin{equation}
\label{f2_86}
d := 1 + b_1 + b_2 + b_3 + b_1 b_2 + b_1 b_3 + b_2 b_3 + b_1 b_2 b_3,\quad \widehat\alpha = 1 + \alpha_1 + \alpha_2.
\end{equation}

Assume that $z\not= 0$.
We can rewrite the differential operator $[...]$ in~(\ref{f2_79}) as a sum of powers of $\theta$, and divide the whole
equality by z to obtain:
$$ \left[
\frac{d}{dz}
(\theta^3 + (b_1 + b_2 + b_3) \theta^2 + (b_1 b_2 + b_1 b_3 + b_2 b_3)\theta + b_1 b_2 b_3
) - \right. $$
\begin{equation}
\label{f2_89}
\left.
- \theta^2 - (\alpha_1 + \alpha_2) \theta - \alpha_1\alpha_2
\right] w(z) = 0,\qquad z\in\mathbb{C}\backslash \{ 0 \}.
\end{equation}
In terms of usual derivatives this relation can be rewritten as
\begin{equation}
\label{f2_90}
z^3 w^{(4)} + c z^2 w''' + (\widehat b - z) z w'' + (d - \widehat\alpha z) w' - \alpha_1\alpha_2 w = 0,\qquad 
   z\in\mathbb{C}\backslash \{ 0 \}. 
\end{equation}
Denote the left-hand side of~(\ref{f2_90}) by $l(z)$. 
Since $w(z)$ is an entire function, then $l(z)$ is entire as well. By continuity we conclude that relation~(\ref{f2_90})
holds for $z=0$.
Set
\begin{equation}
\label{f2_92}
\varphi(z) = \varphi(z;x) :=
\sum_{n=0}^\infty \mathbf{L}_n(x) \frac{ (-1)^n }{ x^n }\frac{z^n}{ n! },\quad z\in\mathbb{C}.
\end{equation}
Then
$$ w(z) = e^{ \frac{z}{x} } \varphi(z),\quad z\in\mathbb{C}. $$
We can calculate the derivatives of $w$ by the Leibniz rule and substitute the resulting expressions
into relation~(\ref{f2_90}). 
If we cancel the term $e^{ \frac{z}{x} }$,
we shall get the following relation:
$$ z^3 \varphi^{(4)} + \frac{4}{x} z^3 \varphi''' + \frac{6}{x^2} z^3 \varphi'' + \frac{4}{x^3} z^3 \varphi' + \frac{1}{x^4} z^3 \varphi + $$
$$ + c z^2 \varphi''' + c \frac{3}{x} z^2 \varphi'' + c \frac{3}{x^2} z^2 \varphi' + c \frac{1}{x^3} z^2 \varphi + $$
$$ + (\widehat b - z) z \varphi'' + (\widehat b - z) \frac{2}{x} z \varphi' + (\widehat b - z) \frac{1}{x^2} z \varphi + $$
\begin{equation}
\label{f2_94}
+ (d - \widehat\alpha z) \varphi' + (d - \widehat\alpha z) \frac{1}{x} \varphi - \alpha_1\alpha_2 \varphi = 0,\qquad 
z\in\mathbb{C}.  
\end{equation}
Denote the left-hand side of~(\ref{f2_94}) by $\widehat l(z)$. 
Observe that
$$ \varphi'(z) = 
\sum_{n=0}^\infty \frac{ (-1)^{n+1} }{ n! } \frac{ \mathbf{L}_{n+1}(x) }{ x^{n+1} }  z^n,\quad  
\varphi''(z) = 
\sum_{n=0}^\infty \frac{ (-1)^{n} }{ n! } \frac{ \mathbf{L}_{n+2}(x) }{ x^{n+2} }  z^n,  $$
$$ \varphi'''(z) = 
\sum_{n=0}^\infty \frac{ (-1)^{n+1} }{ n! } \frac{ \mathbf{L}_{n+3}(x) }{ x^{n+3} }  z^n,\quad
\varphi^{(4)}(z) = 
\sum_{n=0}^\infty \frac{ (-1)^{n} }{ n! } \frac{ \mathbf{L}_{n+4}(x) }{ x^{n+4} }  z^n.  $$
We can substitute the latter expressions into relation~(\ref{f2_94}) to get a series expansion of $\widehat l(z)$, which is equal to zero.
Thus, every Taylor coefficient $\widehat l_k$ is zero, and this provides a recurrence relation for polynomials $\mathbf{L}_n$.

\begin{theorem}
\label{t2_4}
Let $\alpha_1,\alpha_2,\beta_1,\beta_2,\beta_3\in(0,+\infty)$. Consider polynomials
$$ \mathbf{L}_n(x) = \mathbf{L}_n(x;\alpha_1,\alpha_2;\beta_1,\beta_2,\beta_3) =
 {}_{3} F_{3} (-n,\alpha_1,\alpha_2;\beta_1,\beta_2,\beta_3;x),\quad n\in\mathbb{Z}_+, $$
with $\mathbf{L}_{-1}(x) = \mathbf{L}_{-2}(x) = \mathbf{L}_{-3}(x) \equiv 0$.
Let $b_1, b_2, b_3, c, \widehat b, d, \widehat\alpha$ be defined as in~(\ref{f2_82})-(\ref{f2_86}).
The following five-term recurrence relation holds:
$$ \left(
-k(k-1)(k-2) - k(k-1) c - k \widehat b - d
\right)
\mathbf{L}_{k+1}(x) + $$
$$ + \left(
4k(k-1)(k-2) + 3k(k-1) c + 2 k \widehat b + d
\right)
\mathbf{L}_{k}(x) + $$
$$ + \left(
-6k(k-1)(k-2) - 3k(k-1) c - k \widehat b
\right)
\mathbf{L}_{k-1}(x) + $$
$$ + \left(
4k(k-1)(k-2) + k(k-1) c 
\right)
\mathbf{L}_{k-2}(x) -
k(k-1)(k-2) \mathbf{L}_{k-3}(x) = $$
$$ = x
\left[
( k(k-1) + k\widehat\alpha + \alpha_1\alpha_2 )\mathbf{L}_{k}(x) - \right. $$
\begin{equation}
\label{f2_95}
- \left. ( 2k(k-1) + k\widehat\alpha ) \mathbf{L}_{k-1}(x)
+ k(k-1) \mathbf{L}_{k-2}(x)
\right],\qquad   k\in\mathbb{Z}_+. 
\end{equation}

\end{theorem}

\noindent
\textbf{Proof.} Calculate the Taylor coefficients $\widehat l_k$ of $\widehat l(z)$, as it was explained before the statement of the theorem.
Then multiply $\widehat l_k$ by $(-1)^k k! x^{k+1}$ to get relation~(\ref{f2_95}). 
$\Box$

In conditions of Theorem~\ref{t2_4} we additionally assume that
\begin{equation}
\label{f2_97}
\beta_1,\beta_2,\beta_3 \in [1,+\infty). 
\end{equation}
Then
parameters $b_1,b_2,b_3; c,\widehat b, d$ are positive.
This fact ensures that the coefficient by $\mathbf{L}_{k+1}(x)$ in the recurrence relation~(\ref{f2_95}) is non-zero for $k\geq 3$.
Since the coefficient by $\mathbf{L}_{k-3}(x)$ is also non-zero for $k\geq 3$, the recurrence relation~(\ref{f2_95})
is non-trivial in this case.

Notice that by~(\ref{f2_37}) we may write
$$ \mathcal{L}_n(x;\alpha,\delta_1,\delta_2,\kappa_1,\kappa_2) = $$
\begin{equation}
\label{f2_97_1}
= \mathbf{L}_n(x;\delta_1 +1,\delta_2 +1;\alpha+1,\kappa_1+\delta_1 +1,\kappa_2+\delta_2 +1),\quad
n\in\mathbb{Z}_+, 
\end{equation}
where
$\alpha,\delta_1,\delta_2\in(-1,+\infty)$, and
$\kappa_1,\kappa_2\in\mathbb{Z}_+$, are arbitrary parameters.
Therefore one can write the above recurrence relation
for $\mathcal{L}_n(x;\alpha,\delta_1,\delta_2,\kappa_1,\kappa_2)$.

Let us turn to the question about the location of zeros of polynomials $\mathbf{P}_n$ and $\mathbf{L}_n$.

\begin{proposition}
\label{p2_3}
Let $p,q\in\mathbb{Z}_+$: $p\geq q+1$, and
$$ a\in(-1,+\infty);\ \alpha_1,\dots,\alpha_p;\beta_1,\dots,\beta_q\in (0,+\infty), $$ 
are some parameters. 
If
\begin{equation} 
\label{f2_107}
\alpha_j \geq \beta_j,\ j\in\mathbb{Z}_{1,q};\quad \alpha_k \geq 1,\ k\in\mathbb{Z}_{q+1,p},
\end{equation}
then all zeros of polynomials $\mathbf{P}_n(x) = \mathbf{P}_n(x;a,\alpha_1,\dots,\alpha_p;\beta_1,\dots,\beta_q)$
and all zeros of polynomials $\mathbf{L}_n(x) = \mathbf{L}_n(x;\alpha_1,\dots,\alpha_p;\beta_1,\dots,\beta_q)$
lie in the unit disc $\mathbb{D}$.

\end{proposition}

\noindent
\textbf{Proof.} Fix an arbitrary $n\in\mathbb{N}$.
Since 
$$ \mathbf{P}_n(x;a,\alpha_1,\dots,\alpha_p;\beta_1,\dots,\beta_q) =
{}_{p+2} F_{q} (-n,n+a,\alpha_1,\dots,\alpha_p;\beta_1,\dots,\beta_q;x) = $$
$$ = \sum_{k=0}^n (-n)_k (n+a)_k \frac{ (\alpha_1)_k\dots (\alpha_p)_k }{ (\beta_1)_k\dots (\beta_q)_k }
\frac{ x^k }{ k! } = $$
$$ = \sum_{k=0}^n \frac{ n! }{ (n-k)! } (n+a)_k \frac{ (\alpha_1)_k\dots (\alpha_p)_k }{ (\beta_1)_k\dots (\beta_q)_k }
\frac{ (-x)^k }{ k! } = \sum_{k=0}^n d_k z^k =: p(z), $$
where
$$ d_k := \frac{ n! }{ (n-k)! } (n+a)_k \frac{ (\alpha_1)_k\dots (\alpha_p)_k }{ (\beta_1)_k\dots (\beta_q)_k }
\frac{1}{k!} > 0,\quad
z := -x. $$
Thus, the polynomial $p(z)$ has degree $n$ and positive coefficients. The reversed polynomial:
$$ p^*(z) := z^n p(1/z), $$
has degree $n$ and positive coefficients as well.
Observe that
$$ d_k/d_{k+1} = \frac{1}{ (n-k) } \frac{1}{ (n+a+k) } \frac{ (\beta_1+k)\dots (\beta_q+k) }{ (\alpha_1+k)\dots (\alpha_p+k) }
(k+1)  \leq 1,\quad k\in\mathbb{Z}_{0,n-1}, $$
where we used condition~(\ref{f2_107}).
We can apply the Enestr\"om--Kakeya Theorem~(\cite[p. 136]{cit_8100_Marden}) for the polynomial $p^*(z)$
to obtain that all its zeros lie in the domain $D_e := \{ z\in\mathbb{C}:\ |z|>1 \}$.
Therefore the zeros of $\mathbf{P}_n$ lie in $\mathbb{D}$.

We may proceed for polynomials $\mathbf{L}_n$ in a similar way:
$$ \mathbf{L}_n(x;\alpha_1,\dots,\alpha_p;\beta_1,\dots,\beta_q) =
{}_{p+1} F_{q} (-n,\alpha_1,\dots,\alpha_p;\beta_1,\dots,\beta_q;x) = $$
$$ = \sum_{k=0}^n (-n)_k \frac{ (\alpha_1)_k\dots (\alpha_p)_k }{ (\beta_1)_k\dots (\beta_q)_k }
\frac{ x^k }{ k! } = $$
$$ = \sum_{k=0}^n \frac{ n! }{ (n-k)! } \frac{ (\alpha_1)_k\dots (\alpha_p)_k }{ (\beta_1)_k\dots (\beta_q)_k }
\frac{ (-x)^k }{ k! } = \sum_{k=0}^n \widehat d_k z^k =: \widehat p(z), $$
where
$$ \widehat d_k := \frac{ n! }{ (n-k)! } \frac{ (\alpha_1)_k\dots (\alpha_p)_k }{ (\beta_1)_k\dots (\beta_q)_k }
\frac{1}{k!} > 0,\quad
z := -x. $$
Since
$$ \widehat d_k/\widehat d_{k+1} \leq 1,\qquad k\in\mathbb{Z}_{0,n-1}, $$
by the Enestr\"om--Kakeya Theorem we conclude that the reversed polynomial $\widehat p^*$ has its zeros in $D_e$.
Thus, the zeros of $\mathbf{L}_n$ lie in $\mathbb{D}$ as well.
$\Box$

\vspace{1.5cm}

V. N. Karazin Kharkiv National University \newline\indent
School of Mathematics and Computer Sciences \newline\indent
Department of Higher Mathematics and Informatics \newline\indent
Svobody Square 4, 61022, Kharkiv, Ukraine

Sergey.M.Zagorodnyuk@gmail.com; zagorodnyuk@karazin.ua

}

\begin{thebibliography}{99}


\bibitem{cit_DG_1986}
J.J. Duistermaat,  F.A. Gr\"unbaum, \emph{Differential equations in the spectral parameter}, Comm. Math. Phys. 103 (1986), no. 2, pp. 177--240.


\bibitem{cit_110_Duran}
A.J. Dur\'an, M.D. de la Iglesia, \emph{Differential equations for discrete Jacobi-Sobolev orthogonal polynomials}, 
J. Spectr. Theory 8 (2018), no. 1, pp. 191--234.


\bibitem{cit_50_Bateman}
A. Erd\'elyi, W. Magnus, F. Oberhettinger,F.G. Tricomi, 
\emph{Higher transcendental functions. Vols. I, II. Based, in part, on notes left by Harry Bateman}, 
McGraw-Hill Book Company, Inc., New York-Toronto-London, 1953.


\bibitem{cit_70_Bateman}
A. Erd\'elyi, W. Magnus, F. Oberhettinger, F.G. Tricomi, 
\emph{Higher transcendental functions. Vol. III. Based, in part, on notes left by Harry Bateman}, McGraw-Hill Book Co., Inc., 
New York-Toronto-London, 1955.

\bibitem{cit_110_Everitt_2001}
W.N. Everitt, K.H. Kwon, L.L. Littlejohn, R. Wellman,  
\emph{Orthogonal polynomial solutions of linear ordinary differential equations}, 
Proceedings of the Fifth International Symposium on Orthogonal Polynomials, Special Functions and their Applications (Patras, 1999), 
J. Comput. Appl. Math. 133 (2001), no. 1-2, pp.~85--109.



\bibitem{cit_90_Fichtenholz}
G.M. Fichtenholz, \emph{Infinite series: ramifications}, Revised English edition. 
Translated from the Russian and freely adapted by Richard A. Silverman. 
The Pocket Mathematical Library, Course 4. Gordon and Breach Science Publishers, New York-London-Paris, 1970.


\bibitem{cit210_Horozov}
E. Horozov, \emph{Vector orthogonal polynomials with Bochner's property}, Constr. Approx. 48 (2018), no. 2, pp. 201--234.

\bibitem{cit_5000_Ismail}
M.E.H. Ismail,  \emph{Classical and quantum orthogonal polynomials in one variable}, With two chapters by Walter Van Assche. 
With a foreword by Richard A. Askey. Encyclopedia of Mathematics and its Applications, 98. Cambridge University Press, Cambridge, 2005.


\bibitem{cit_5105_Koekoek_book}
R. Koekoek, P.A. Lesky, R.F. Swarttouw, \emph{Hypergeometric orthogonal polynomials and their 
$q$-analogues}, With a foreword by Tom H. Koornwinder, Springer Monographs in Mathematics, Springer-Verlag, Berlin, 2010.


\bibitem{cit_5150_M_X}
F. Marcell\'an, Yuan Xu, \emph{On Sobolev orthogonal polynomials}, Expo. Math. 33 (2015), no. 3, 308--352.

\bibitem{cit_8100_Marden}
M. Marden, \emph{Geometry of polynomials}, Second edition. Mathematical Surveys, No. 3, American Mathematical Society, Providence, R.I., 
1966.


\bibitem{cit_5700_McBride__1971}
E.B. McBride, \emph{Obtaining generating functions}, Springer Tracts in Natural Philosophy, 
Vol. 21, Springer-Verlag, New York-Heidelberg, 1971.


\bibitem{cit_5150_Rainville}
E.D. Rainville, \emph{Special functions}, Reprint of 1960 first edition, Chelsea Publishing Co., Bronx, N.Y., 1971.

\bibitem{cit_7100__S_Zh}
V. Spiridonov, A. Zhedanov, 
\emph{Classical biorthogonal rational functions on elliptic grids}, 
C. R. Math. Acad. Sci. Soc. R. Can. 22 (2000), no. 2, pp. 70--76.



\bibitem{cit_50000_Gabor_Szego}
G. Szeg\"o, \emph{Orthogonal polynomials}, Fourth edition. 
American Mathematical Society, Colloquium Publications, Vol. XXIII, American Mathematical Society, Providence, R.I., 1975.

\bibitem{cit_80000_Zagorodnyuk_JAT_2020}
S.M. Zagorodnyuk, \emph{On some classical type Sobolev orthogonal polynomials}, J. Approx. Theory 250 (2020), 105337, 14 pp.

\bibitem{cit_80000_Zagorodnyuk_JDEA_2021}
S.M. Zagorodnyuk, \emph{On some Sobolev spaces with matrix weights and classical type Sobolev orthogonal polynomials},
J. Difference Equ. Appl. 27 (2021), no. 2, pp. 261--283.


\bibitem{cit_90000_Zagorodnyuk_AOT_2022}
S.M. Zagorodnyuk, \emph{On the multiplication operator by an independent variable in matrix Sobolev spaces}, Adv. Oper. Theory 
7 (2022), no. 4, Paper No. 54.



\end{thebibliography}
\end{document}